\documentclass{amsart}
\usepackage{amsmath} 
\usepackage{amssymb}
\usepackage{mathrsfs}

\newtheorem{theorem}{Theorem}[section] 
\newtheorem{claim}[theorem]{Claim}

\newtheorem{conclusion}[theorem]{Conclusion}
\newtheorem{observation}[theorem]{Observation}

\theoremstyle{definition}
\newtheorem{definition}[theorem]{Definition}

\newtheorem{discussion}[theorem]{Discussion}

\theoremstyle{remark}
\newtheorem{remark}[theorem]{Remark}

\newtheorem{notation}[theorem]{Notation}

\newcommand{\CH}{{\rm CH}}

\newcommand{\suc}{{\rm suc}}

\newcommand{\Dom}{{\rm Dom}}
\newcommand{\Rang}{{\rm Rang}}

\newcommand{\rest}{{\restriction}}

\newcommand{\wilog}{{\rm without loss of generality}}

\newcommand{\then}{{\underline{then}}}
\newcommand{\when}{{\underline{when}}}
\newcommand{\where}{{\underline{where}}}
\newcommand{\Then}{{\underline{Then}}}

\newcommand{\Iff}{{\underline{iff}}}
\newcommand{\mn}{{\medskip\noindent}}
\newcommand{\sn}{{\smallskip\noindent}}

\newcommand{\cA}{{\mathscr A}}

\newcommand{\gh}{{\mathfrak h}}

\newcommand{\cH}{{\mathscr H}}

\newcommand{\cJ}{{\mathscr J}}

\newcommand{\cI}{{\mathscr I}}

\newcommand{\cW}{{\mathscr W}}
\newcommand{\bbP}{{\mathbb P}}
\newcommand{\cP}{{\mathscr P}}

\newcommand{\varp}{{\varepsilon}}

\newcommand{\bbQ}{{\mathbb Q}}
\newcommand{\bbR}{{\mathbb R}}

\newcommand{\cT}{{\mathscr T}}
 
\newcommand{\cU}{{\mathscr U}}

\newcommand{\cf}{{\rm cf}}

\newcount\skewfactor
\def\mathunderaccent#1#2 {\let\theaccent#1\skewfactor#2
\mathpalette\putaccentunder}
\def\putaccentunder#1#2{\oalign{$#1#2$\crcr\hidewidth
\vbox to.2ex{\hbox{$#1\skew\skewfactor\theaccent{}$}\vss}\hidewidth}}
\def\name{\mathunderaccent\tilde-3 }

\newenvironment{PROOF}[2][\proofname.]
   {\begin{proof}[#1]\renewcommand{\qedsymbol}{\textsquare$_{\rm #2}$}}
   {\end{proof}}

\begin{document}

\title {Preserving old $([\omega]^{\aleph_0},\supseteq^*)$ is proper}
\author {Saharon Shelah}
\address{Einstein Institute of Mathematics\\
Edmond J. Safra Campus, Givat Ram\\
The Hebrew University of Jerusalem\\
Jerusalem, 91904, Israel\\
 and \\
 Department of Mathematics\\
 Hill Center - Busch Campus \\ 
 Rutgers, The State University of New Jersey \\
 110 Frelinghuysen Road \\
 Piscataway, NJ 08854-8019 USA}
\email{shelah@math.huji.ac.il}
\urladdr{http://shelah.logic.at}
\thanks{Research supported by the United States-Israel Binational
Science Foundation (Grants No. 2002323 and 2006108) and the NSF.
The author thanks Alice Leonhardt for the beautiful typing.  First
 typed Dec. 19, 2007. Paper 960 on author's list}

 %formerly f887 - first typed 07/Dec/19

 %previous version - 2017/Sept/8



\date{September 18, 2017}

\begin{abstract}
We give some sufficient and necessary conditions on a forcing notion
$\bbQ$ for preserving the forcing notion
$([\omega]^{\aleph_0},\supseteq^*)$ being proper.  They cover many
reasonable forcing notions.
\end{abstract}

\maketitle
\numberwithin{equation}{section}
\setcounter{section}{-1}
\newpage

\section* {Anotated Content}
\bigskip

\noindent
\S0 \quad Introduction, pg.\pageref{intro}
\mn
\begin{enumerate}
\item[${{}}$]   [I.e. Definition \ref{1a.1}, we define the problem
and some variants.]
\end{enumerate}
\bigskip

\noindent
\S1 \quad Properness of $\bbP_{{\cA}[\bold V]}$ and CH, pg.\pageref{prop}
\mn
\begin{enumerate}
\item[${{}}$]   [Under CH, if non-meagerness of $({}^\omega
2)^{\bold V}$ is preserved then $\bbP_{{\cA}_*[\bold V]}$ is proper,
(\ref{2b.1}).  If $\bold V$ fails to satisfy CH, then usually 
$\bbP_{{\cA}_*[\bold V]}$ is not proper after a forcing adding a 
new real and satisfying a relative of being proper, e.g. satisfies
c.c.c. or is any true creature forcing.]
\end{enumerate}
\bigskip

\noindent
\S2 \quad General sufficient conditions, pg. \pageref{general}
\mn
\begin{enumerate}
\item[${{}}$]  [If $\bold V$ satisfies CH and $\bbQ$ is
 c.c.c. then $\Vdash_{\bbQ} ``{\bbP}_{{\cA}[\bold V]}$ is proper", see
 in \ref{3c.1}.  In \ref{3c.3} we replace ${\cA}^{\bold V}_*$ by a 
forcing notion
$\bbR$ adding no $\omega$-sequence, $\bbQ$ is c.c.c. even in
$\bold V^{\bbP}$.  Instead ``$\bbQ$ satisfies the c.c.c." it
  suffices to demand $\bbQ$ satisfies a weaker condition.  
Lastly, in \ref{3c.7} we prove some proper forcing does not preserve.]
\end{enumerate}
\newpage

\section {Introduction} \label{intro}

We investigate the question ``$\Pr^+_1(\bbQ,\bbR)$", which means 
that the proper forcing $\bbQ$ preserves that the (old) ${\bbR}$ is proper for
various $\bbR$'s.  In what follows, $B \subseteq^* A$ means $|B
\backslash A| < \aleph_0$, and $A \supseteq^* B$ means the same.

Recall:

\begin{definition}
\label{1a.0}
properness:
\mn
\begin{enumerate}
\item[(a)]  Assume that $N \prec (\cH(\chi),\in),\bbP \in N$ is a
  forcing notion and $q \in \bbP$.  We say that $q$ is
  $(N,\bbP)$-generic \Iff \, for every dense $D \subseteq \bbP$, if $D
  \in N$ then $D \cap N$ is pre-dense above $q$.
\sn
\item[(b)]  A forcing notion $\bbP$ is proper \Iff \, for every
  sufficiently large regular $\chi$ and every countable $N \prec
  (\cH(\chi),\in)$, if $p,\bbP \in N$ then there is a condition $q \in
  \bbP,q \ge p$ such that $q$ is $(N,\bbP)$-generic.
\end{enumerate}
\mn
Gitman proved that $\Pr^+_1(\bbQ,\bbP_{{\cP}(\omega)[\bold V]})$
(see definition below, \where \, $\bbP_{\cP(\omega)^{[\bold V]}}$ 
is the forcing notion $(\{A \in \bold V:A \subseteq \omega,|A| = 
\aleph_0\},\supseteq^*)$, \when \, $\bbQ$ is adding Cohen reals 
(or just Cohen subsets
even $> 2^{\aleph_0}$ many).  But no other examples were
known even Sacks forcing.  Also for e.g. $\bold V \models ``V = L"$, we
did not know a forcing making it not proper.

We thank Victoria Gitman for asking us the question and Otmar Spinas
and Haim Horowitz for comments and Shimoni Garti for many more.

Let us state the problem and relatives.  We are interested mainly in
the case $\bbQ$ is proper.
\end{definition}

\begin{definition}
\label{1a.1} 
1) Let Pr$_1(\bbQ,\bbP)$ means: $\bbQ,\bbP$ are forcing notions 
and $\Vdash_{\bbQ} ``\bbP$, i.e. 
$\bbP^{\bold V}$ is a proper forcing".

\noindent
1A) Let $\Pr^+_1(\bbQ,\bbP)$ be defined similarly but adding ``$\bbQ$
is proper".

\noindent
2) For ${\cA} \subseteq {\cP}(\omega)$ let $\bbP_{\cA}$ be
${\cA} \backslash [\omega]^{< \aleph_0}$ ordered by $\supseteq^*$, inverse
almost inclusion.

\noindent
3) Let ${\cA}_* = \cA_*[\bold V] = ([\omega]^{\aleph_0})^{\bold V}$.
\end{definition}

\begin{observation}
\label{1a.3}  A necessary condition for $\Pr_1(\bbQ,\bbP)$ is: 
\mn
\begin{enumerate}
\item[$(*)_1$]   if $\chi$ is regular and
large enough, $N \prec ({\cH}(\chi),\in)$ 
is countable, $\bbQ,\bbP \in N,q_1 \in \bbQ$ is
$(N,\bbQ)$-generic and $r_1 \in N \cap \bbP$ \then \, we can find
$(q_2,r_2)$ such that:
\sn
\begin{enumerate}
\item[$\odot$]  
\begin{enumerate}
\item[(a)]  $q_1 \le_{\bbQ} q_2$
\sn
\item[(b)]  $r_1 \le_{\bbP} r_2$
\sn
\item[(c)]  $q_2 \Vdash ``r_2$ is $(N[\name G_{\bbQ}],\bbP)$-generic".
\end{enumerate}
\end{enumerate}
\end{enumerate}
\end{observation}

\begin{definition}
\label{1a.7}  
1) We define $\Pr^-(\bbQ,\bbP) = \Pr_2(\bbQ,\bbP)$ as the 
necessary condition from \ref{1a.3}.

\noindent
2) Let $\Pr_3(\bbQ,\bbP)$ mean that $\bbQ,\bbP$ are
forcing notions and for some $\lambda$ and stationary $S \subseteq
[\lambda]^{\aleph_0}$ from $\bold V$ we have $\Vdash_{\bbQ} 
``\Bbb P$ is $S$-proper", and note that $S$ remains stationary of course.

\noindent
3) Pr$_4(\bbQ,\bbP)$ is defined similarly but $S \in \bold V^{\bbQ}$, still $S
   \subseteq ([\lambda]^{\aleph_0})^{\bold V}$, so $S$ is actually
   $\name S$, a $\bbQ$-name.

\noindent
4) Pr$_5(\bbQ,\bbP)$ is the statement (A) of \ref{1a.10}(4) below.

\noindent
5) Let $\Pr^+_\ell(\bbQ,\bbP)$ means $\Pr_\ell(\bbQ,\bbP)$ and $\bbQ$
   is a proper forcing, for $\ell=2,3,4,5$.
\end{definition}

\begin{claim}
\label{1a.10}  
1) $\Pr_2(\bbQ,\bbP)$ means that for $\lambda$ large enough,
 letting $S = ([\lambda]^{\aleph_0})^{\bold V}$, we have
$\Vdash_{\bbQ} ``\bbP$ is $S$-proper".

\noindent
2) {\rm Pr}$_1(\bbQ,\bbP) \Rightarrow \text{\rm Pr}_2(\bbQ,\bbP) 
\Rightarrow \text{\rm Pr}_3(\bbQ,\bbP)$; similarly for $\Pr^+$.

\noindent
3) Also {\rm Pr}$_3(\bbQ,\bbP) \Rightarrow \text{\rm Pr}_4(\bbQ,\bbP) 
\Rightarrow \text{\rm Pr}_5(\bbQ,\bbP)$; similarly for $\Pr^+$.

\noindent
4) If $\bbQ,\bbP$ are forcing notions, $\chi$ large enough and
regular, \then \, (A) $\Leftrightarrow$ (B) where
\mn
\begin{enumerate}
\item[$(A)$]   for some countable $N \prec ({\cH}(\chi),\in)$ and for
   some $q \in \bbQ,p \in \bbP$ we have
\sn
\begin{enumerate}
\item[$(a)$]   $q$ is $(N,\bbQ)$-generic
\sn
\item[$(b)$]   $q \Vdash_{\bbQ} ``p$ is $(N[\name G_{\bbQ}],\bbP)$-generic"
\end{enumerate}
\sn
\item[$(B)$]   for some $q_* \in \bbQ,p_* \in \bbP$ we have 
$\Pr_4(\bbQ_{\ge q_*},\bbP_{\ge p_*})$.
\end{enumerate}
\end{claim}

\begin{PROOF}{\ref{1a.10}}
 Easy.  
\end{PROOF}

\begin{notation}
\label{z2}
$<^*_\chi$ denotes a well ordering of $\cH(\chi)$.
\end{notation}

\noindent
Recall (Balcar-Pelant-Simon \cite{BPS}, or see, e.g. Blass \cite{Bls10})
\begin{definition}
\label{z5}
$\gh$ is the following cardinal invariant, it is the minimal
cardinality $\chi$ (necessarily regular) such that forcing with
$\bbP_{\cA_*}$ adds a new sequence of ordinals of length $\chi$.
\end{definition}

\begin{notation}
\label{z8}
If $\cT$ is a tree, then $\suc_{\cT}(p)$ is the set of immediate
successors of $p \in \cT$ in the tree order.
\end{notation}
\newpage

\section {Properness of $\bbP_{{\cA}_*[\bold V]}$ and CH} \label{prop}

\begin{claim}
\label{2b.1}  
Assume $\bold V_0 \models$ {\rm CH}, $\bold V_1 \supseteq 
\bold V_0$, e.g. $\bold V_1 = \bold V^{\bbQ}_0$ and
let ${\cA} = {\cA}_*[\bold V_0]$.
\mn
\begin{enumerate}
\item[(a)]  If $\aleph^{\bold V_0}_1$ is a countable ordinal in $\bold
  V_1$, then $\bold V_1 \models ``\bbP_{\cA}$ is proper".
\sn
\item[(b)]  If $\aleph^{\bold V_0}_1 = \aleph^{\bold V_1}_1$ and $\bold
  V_1 \models ``({}^\omega 2)^{\bold V_0}$ is non-meagre", then $\bold
  V_1 \models ``\bbP_{\cA}$ is proper".
\end{enumerate}
\mn
In both cases, if $\bold V_1$ is a generic extension of $\bold V_0$ by
the forcing notion $\bbQ$ then it means that $\Pr_1(\bbQ,\bbP_{\cA})$ holds.
\end{claim}

\begin{PROOF}{\ref{2b.1}}
Assume that $\bold V_1 \supseteq \bold V_0$.

If $\bold V_1 \models ``\aleph^{\bold V_0}_1$ is
countable" then recalling $\bold V_0 \models \text{ CH}$ clearly
$\bold V_1 \models ``{\cA}$ is countable" so we
know that $\bbP_{\cA}$ is proper in $\bold V_1$, thus proving clause (a).  
So from now on we assume $\aleph^{\bold V_0}_1$ is not collapsed.

In $\bold V_0$ let $\cT = {}^{\omega_1>}(\omega_1)$ and choose a
subset $\cA' \subseteq \cA$ such that $\cA'$ is $\subseteq^*$-dense in
$\cA$ and $(\cA',\supseteq^*)$ is tree-isomorphic to $\cT$.  Let $\pi$ be
the isomorphism between these trees\footnote{this is trivial 
as $\bold V_0 \models$ CH, however always there is a dense tree with
${\frak h}$ levels by the celebrated theorem of Balcar-Pelant-Simon}.
Notice that all this is done in $\bold V_0$ (recalling that $\bold V_0
\models \CH$).
In $\bold V_0$ there is a
sequence $\bar{\cT} = \langle {\cT}_\alpha:\alpha <
\omega_1\rangle$ which is $\subseteq$-increasing continuous with union
${\cT}$ and each ${\cT}_\alpha$ countable.  Also there is $\bar
C = \langle C_\delta:\delta < \omega_1,\delta$ is a limit 
ordinal$\rangle \in \bold V_0$ such that $C_\delta \subseteq \delta =
\sup(C_\delta)$, otp$(C_\delta) = \omega$.  Let ${\cT}'_\delta = 
{\cT}_\delta \rest \{\eta \in {\cT}_\delta:\ell g(\eta) \in C_\delta\}$. 

In $\bold V_1$ choose a sufficiently large regular cardinal $\chi$, and
 let $N \prec ({\cH}(\chi),\in)$ be countable such that $\cA,\pi,
\bar{\cT} \in N$ and let $\delta = \omega_1 \cap N$, clearly
${\cT} \cap N = {\cT}_\delta$.  We have to prove the statement:
\mn
\begin{enumerate}
\item[$(*)_0$]   ``for every $p \in \bbP_{\cA} \cap N$ there is 
$q \in \bbP_{\cA}$ above $p$ which is $(N,\bbP_{\cA})$-generic". 
\end{enumerate}
\mn
As $\bold V_0 \models \CH$ and the density of $\cA'$ in $\cA$ and
$(\cA',\supseteq^*)$ being isomorphic in $\bold V_0$ by $\pi$ to 
${\cT}$ this is equivalent (in $\bold V_1$, of course) to:
\mn
\begin{enumerate}
\item[$(*)_1$]   for every $\nu \in {\cT} \cap N = 
{\cT}_\delta$ there is $\eta \in {\cT}$ which is 
$(N,{\cT})$-generic and $\nu \le_{\cT} \eta$.
\end{enumerate}
\mn
In $\bold V_0$ we let $\bar S = \langle S_\delta:\delta < \omega_1$ a limit
ordinal$\rangle$ where $S_\delta = \{\bar \nu:\bar \nu = \langle \nu_n:n <
\omega\rangle$ is $<_{\cT}$-increasing, $\nu_n \in {\cT}'_\delta$, 
moreover $\ell g(\nu_n)$ is the $n$-th member of $C_\delta\}$.

As $(\forall \nu \in {\cT}_\delta)(\exists \rho)(\nu <_{\cT}
\rho \in {\cT}'_\delta)$, and [$\bar\nu \in S_\delta \Rightarrow$
  there is a $<_{\cT}$-upper bound $\rho \in \cT$ of $\bar\nu$, in $\bold
  V_0$, of course] recalling $\cT_\delta,S_\delta \in 
\bold V_0$ clearly $(*)_1$ is equivalent (in $\bold V_1$, of course) to
\mn
\begin{enumerate}
\item[$(*)_2$]   for every $\nu \in {\cT}'_\delta$ there is $\bar
\nu \in S_\delta$ such that $\nu \in \text{ Rang}(\bar\nu)$ and
$\bar \nu$ induce a subset of ${\cT}_\delta$ generic over $N$ (i.e.
$(\forall A)[A \in N \text{ is a dense open subset of } {\cT}
\Rightarrow A \cap \{\nu_n:n < \omega\} \ne \emptyset]$.
\end{enumerate}
\mn
Now a sufficient condition for $(*)_2$ is
\mn
\begin{enumerate}
\item[$(*)_3$]   $S_\delta$, as a set of $\omega$-branches of the
tree ${\cT}'_\delta$, is non-meagre.
\end{enumerate}
\mn
But in $\bold V_0,{\cT}'_\delta$ and ${}^{\omega >}\omega$ are
isomorphic and $S_\delta$ is the set of all $\omega$-branches of
${\cT}'_\delta$, so by an assumption from part (b), 
$(*)_3$ holds so we are done.
\end{PROOF}

\begin{discussion}
\label{2b.4}  
However, there can be ${\cA} \subseteq {\cP}(\omega)$ such 
that $({\cA},\subseteq^*)$ is a variation of Souslin tree.  
\end{discussion}

\begin{claim}
\label{2b.7}  
1) We have {\rm Pr}$_1(\bbQ,\bbP_{\cA_*[\bold V]})$ \when \,:
\mn
\begin{enumerate}
\item[$(a)$]   $\aleph_1^{\bold V[\bbQ]} = \aleph_1$
\sn
\item[$(b)$]   $\Vdash_{\bbQ} ``|\lambda| = \aleph_1$ where
$\lambda = (2^{\aleph_0})^{\bold V}"$
\sn
\item[$(c)$]   moreover letting $\langle \name u_i:i
< \aleph_1\rangle$ be a $\bbQ$-name of a
$\subseteq$-increasing continuous sequence of countable
subsets of $\lambda$ with union $\lambda$, the $\bbQ$-name $\name S = 
\{i:u_i \in \bold V\}$ is forced to contain a club (of $\aleph_1$)
\sn
\item[$(d)$]   forcing with $\bbQ$ preserves 
``$({}^\omega 2)^{\bold V}$ is non-meagre".
\end{enumerate}
\mn
2) Assume the forcing notion $\bbQ$ satisfies (a) + (d), 
{\rm Pr}$_4(\bbQ,\bbP_{{\cA}_*[\bold V]})$ as witnessed by $S$ and 
$\bbQ$ is proper and $\name S$ is forced to be stationary.

\Then \, the forcing notion 
$\bbQ * \text{\rm Levy}(\aleph_1,(|\bbQ|^{\aleph_0})^{\bold V}) * 
\bbQ_{\name S}$ preserves ``$\bbP_{{\cA}_*[\bold V]}$ is proper"
where $\bbQ_S$ is the (well known) shooting of a club through the 
stationary subsets of $\omega_1$ (to make clause (c) hold).
\end{claim}

\begin{PROOF}{\ref{2b.7}}
Like \ref{2b.1}. 
\end{PROOF}

\noindent
In what follows we prove that many forcing notions destroy
properness. We need a preliminary concept.
\begin{definition}
\label{2b.13} 
For $\lambda > \kappa$ we say that a forcing notion
$\bbQ$ is $(\lambda,\kappa)$-newly proper (omitting $\kappa$
means $\kappa = \aleph_0$ and we define $(\lambda,<\chi)$-newly proper 
similarly) when: if $\bar N =
\langle (N_\eta,\nu_\eta):\eta \in {}^{\omega >}\lambda\rangle$
satisfies $\circledast$ below and $\bbQ \in N_{<>}$ and $p \in \bbQ \cap
N_{<>}$ \then \, we can find $q,\name \eta$ such that $\boxtimes$ below
holds where:
\mn
\begin{enumerate}
\item[$\circledast$]   for some cardinal $\chi > \lambda$
\sn
\begin{enumerate}
\item[$(a)$]    $N_\eta \prec ({\cH}(\chi),\in,<^*_\chi)$ is countable
\sn
\item[$(b)$]   if $\nu \triangleleft \eta$ then $N_\nu \prec N_\eta$
\sn
\item[$(c)$]   $N_{\eta_1} \cap N_{\eta_2} = N_{\eta_1 \cap \eta_2}$
if $\kappa = \aleph_0$ and $N^\kappa_{\eta_1} \cap N^\kappa_{\eta_2} = 
N^\kappa_{\eta_1 \cap \eta_2}$ generally where $N^\kappa_\eta :=
\cup\{v \in N_\eta:|v| \le \kappa\}$
\sn
\item[$(d)$]   $\nu_\eta \in N_\eta \backslash \cup\{N^\kappa_{\eta \rest
m}:m < \ell g(\eta)\}$ hence $\nu_\eta \notin
\cup\{N_\nu:\neg(\eta \trianglelefteq \nu)$ and $\nu \in 
{}^{\omega >}\lambda\}$ 
\sn
\item[$(e)$]   $\nu_\eta \in {}^{\ell g(\eta)}\lambda$ and $\ell
< \ell g(\eta) \Rightarrow \nu_{\eta \rest \ell} \trianglelefteq \nu_\eta$
\end{enumerate}
\sn
\item[$\boxtimes$]
\begin{enumerate}
\item[(a)]  $p \le_{\bbQ} q$
\sn
\item[(b)]  $q \Vdash_{\bbQ} ``\cup\{N_{\name \eta \rest n}
[\name {\bold G}_{\bbQ}]:n < \omega\} \cap \bold V =
\cup\{N_{\name \eta \rest n}:n < \omega\}"$
\sn
\item[(c)]  $q \Vdash_{\bbQ} 
``\name \eta \in {}^\omega \lambda$ is new, i.e. $\name \eta \notin
({}^\omega \lambda)^{\bold V}"$
\sn
\item[(c)$^+$]   moreover if $\kappa > \aleph_0$ and $\cT
\in \bold V$ is a sub-tree of ${}^{\omega >} \lambda$ of cardinality
$\le \kappa$ then $\name \eta \notin \text{\rm lim}(\cT)$,
i.e. $\{\name \eta \rest n:n < \omega\} \notin \cT$.
\end{enumerate}
\end{enumerate}
\end{definition}

\begin{observation}
\label{2b.14}
If $\langle N_\eta:\eta \in {}^{\omega >}\lambda\rangle$ satisfies
clauses (a),(b),(c) of $\circledast$ of Definition \ref{2b.13}, \then
\, the following conditions are equivalent:
\mn
\begin{enumerate}
\item[$\bullet_1$]  there is $\langle \nu_\eta:\eta \in {}^{\omega
    >}\lambda\rangle$ such that clauses (d),(e) of $\circledast$ of
  Definition \ref{2b.13}
\sn
\item[$\bullet_2$]  if $\eta \in {}^{\omega >}\lambda$, then $N_\eta
  \cap \lambda \nsubseteq \cup\{N_{\eta \rest \ell}:\ell < \ell
  g(\eta)\}$.
\end{enumerate}
\end{observation}

\noindent
For a proper forcing notion adding a new real it is quite easy to be
$\aleph_1$-newly proper; e.g.
\begin{claim}
\label{2b.19} 
Assuming $2^{\aleph_0} \ge \lambda = \text{\rm cf}(\lambda) > \aleph_1$,  
sufficient conditions for ``$\bbQ$ is $\lambda$-newly proper" are:
\mn
\begin{enumerate}
\item[$(a)$]   $\bbQ$ is c.c.c. and adds a new real
\sn
\item[$(b)$]   $\bbQ$ is Sacks forcing
\sn
\item[$(c)$]   $\bbQ$ is a tree-like creature forcing in the sense of
Roslanowski-Shelah \cite{RoSh:470}.
\end{enumerate}
\end{claim}

\begin{PROOF}{\ref{2b.19}}
Easy; for clause (a) we use $q=p$ for $\boxplus$ of the definition
noting that: if $\eta \in {}^{\omega >}\lambda$ then $p$ is
$(N_\eta,\bbQ)$-generic. 
For clauses (b),(c) we use fusion but in the $n$-th step use members of
$N_\eta \cap \bbQ$ for $\eta \in {}^n \lambda$, we get as many distinct
$\eta$'s as we can.
\end{PROOF}

\begin{theorem}
\label{2b.10}  
We have
$\Vdash_{\bbQ} ``\bbP_{{\cA}_*[\bold V]} \text{ is not proper}"$ \when :
\mn
\begin{enumerate}
\item[$(a)$]    $\bold V \models 2^{\aleph_0} \ge \aleph_2$
\sn
\item[$(b)$]    $\lambda$ is regular, $\aleph_2 \le \lambda \le
2^{\aleph_0}$ and\footnote{If $\lambda = \aleph_2$ the rest of 
clause (b) follows.}
 $\alpha < \lambda \Rightarrow 
\text{\rm cf}([\alpha]^{\aleph_0},\subseteq) < \lambda$ hence (by
\cite{Sh:420}) there is a stationary ${\cU}_\alpha \subseteq 
[\alpha]^{\aleph_0}$ of cardinality $< \lambda$
\sn
\item[$(c)$]   ${\gh} < \lambda$
\sn
\item[$(d)$]   the forcing notion $\bbQ$ adds at least one real and is 
$\lambda$-newly proper.
\end{enumerate}
\end{theorem}

\begin{PROOF}{\ref{2b.10}}
Let $\chi$ be large enough and for transparency, $x \in {\cH}(\chi)$.  
 
By Rubin-Shelah \cite{RuSh:117}, see more \cite[Ch.XI]{Sh:f} 
in $\bold V$ there is a sequence 
$\langle N_\eta:\eta \in {}^{\omega >}\lambda\rangle$ such that:  
\mn
\begin{enumerate}
\item[$\boxdot_1$] 
\begin{enumerate}
\item[(a)]  $N_\eta \prec ({\cH}(\chi),\in)$
\sn
\item[(b)]  $\bbQ,x \in N_\eta$
\sn
\item[(c)]  $N_\eta$ is countable
\sn
\item[(d)]  $N_{\eta_1} \cap N_{\eta_2} = N_{\eta_1 \cap \eta_2}$.
\end{enumerate}
\end{enumerate}
\mn
Now for each $\eta \in {}^\omega \lambda$ let $N_\eta =
\cup\{N_{\eta \rest k}:k < \omega\}$; we can easily add:
\mn
\begin{enumerate}
\item[(e)]   there is $\cW$ such that:
\sn
\begin{enumerate}
\item[$(\alpha)$]  $\cW$ is a subtree of ${}^{\omega >}\lambda$ 
\sn
\item[$(\beta)$]  $\langle \rangle \in \cW$
\sn
\item[$(\gamma)$]  if $\eta \in \cW$ then $(\exists^\lambda
  \alpha)(\eta \char 94 \langle \alpha \rangle \in \cW)$
\sn
\item[$(\delta)$]  if $\eta \in \lim(W)$ then
$\eta \in {}^\omega \lambda$ is increasing, and
$\sup(N_\eta \cap \lambda) = \sup(\Rang(\eta))$
\sn
\item[$(\varp)$]  we can choose $\nu_\eta \in N_\eta$ for $\nu \in
  \cW$ as in clauses (d),(e) of $\circledast$ of \ref{2b.13}.
\end{enumerate}
\end{enumerate}
\mn
By Balcar-Pelant-Simon \cite{BPS} there is
${\cT} \subseteq [\omega]^{\aleph_0}$ such that
\mn
\begin{enumerate}
\item[$\boxdot_2$]  
\begin{enumerate}
\item[$(\alpha)$]  $(\cT,\supseteq^*)$ is
 a tree with ${\gh}$ levels (${\gh}$ is the cardinal invariant from
 \ref{z5}, a regular cardinal $\in [\aleph_1,2^{\aleph_0}]$), the tree
$\cT$ has a root and each node has $2^{\aleph_0}$ many immediate successors, 
i.e. ${\cT}$ has splitting to $2^{\aleph_0}$)
\sn
\item[$(\beta)$]  ${\cT}$ is dense in
$([\omega]^{\aleph_0},\supseteq^*)$, i.e.  in
$\bbP_{\cP(\omega)^{[\bold V]}} = \bbP_{{\cA}_*[\bold V]}$
recalling \ref{1a.1}(2).
\end{enumerate}
\end{enumerate}
\mn
Choose $\bar h$ such that
\mn
\begin{enumerate}
\item[$\boxdot_3$]   $\bar h = \langle h_p:p \in {\cT}\rangle$
 satisfies:  $h_p$ is a one-to-one function from $\suc_{\cT}(p)$ onto 
$2^{\aleph_0} \backslash \{h_{p_0}(p_1):p_0 <_{\cT} p_1 <_{\cT} p$
and $p_1 \in \suc_{\cT}(p_0)\}$.
\end{enumerate}
\mn
So \wilog \,
\mn
\begin{enumerate}
\item[$\boxdot_4$]   ${\cT} \in N_{<>},\gh \in N_{<>}$ and $\bar h \in N_{<>}$.
\end{enumerate}
\mn
As $\bbQ$ is $\lambda$-newly proper there are $\name \eta,q$ as in
$\boxtimes$ of Definition \ref{2b.13}.
Let $\bold G \subseteq \bbQ$ be generic over $\bold V$ such that
$q \in \bold G$, let $\eta = \name\eta[G]$ and
$M_2 := N_{\name\eta[G]} := \cup\{N_{\eta \rest n}[\bold G]:n
< \omega\}$, so $M_2 \prec ({\cH}(\chi)^{\bold V[\bold G]},
{\cH}(\chi)^{\bold V},\in)$ is countable, pedantically
$(|M_2|,\cH(\chi)^{\bold V} \cap |M_2|, \in \rest |M_2|) 
\prec (\cH(\chi)^{\bold V[\bold G]},
\cH(\chi)^{\bold V},\in \rest \cH(\chi)^{\bold V[\bold G]})$.

By $\boxtimes$ of \ref{2b.13}, i.e. the choice of $\name\eta,q$ as 
$q \in \bold G$ we have $M_1 = M_2 \cap 
{\cH}(\chi)^{\bold V}$ is $\cup\{N_{\eta \rest n}:n < \omega\}$, and of
course $M_1 \prec ({\cH}(\chi),\in)$.  Toward contradiction assume
$\bold V[\bold G] \models ``{\cP}_{{\cA}_*[\bold V]}$ is proper",
hence some $p_* \in \bbP_{{\cA}_*[\bold V]}$ is
$(M_2,\bbP_{\cA_*[\bold V]})$-generic.  But $\cT$ is dense in
$\bbP_{\cA_*[\bold V]}$ so \wilog \, $p_* \in \cT$ and $p_*$
is $(M_2,{\cT})$-generic.

Since $\gh \in N_{<>}$ and $\gh < \lambda$, \wilog \, 
$\eta \in {}^{\omega >}\lambda \Rightarrow N_\eta \cap {\gh} 
= N_{<>} \cap {\gh}$.
For any $\alpha < \lambda$ let 

\[
{\cI}_\alpha = \{p \in {\cT}:\text{ for some }
p_0 \in {\cT} \text{ we have } p \in \text{ suc}_{\cT}(p_0) 
\text{ and } h_{p_0}(p) = \alpha\}
\]

\mn
and letting $\cT_\alpha$ be the $\alpha$-th level of $\cT$ and let

\[
{\cI}^+_\alpha = \{p \in \bbP_{{\cA}_*[\bold V]}:p \text{ is
above some member of } {\cT}_\alpha\}.
\]

\mn
Now clearly (in $\bold V$ and in $\bold V[\bold G]$):
\mn
\begin{enumerate}
\item[$(*)_1$]   $(a) \quad {\cI}_\alpha$ is a pre-dense subset of
${\cT}$ (and of $\bbP_{{\cA}_*[\bold V]}$)
\sn
\item[${{}}$]   $(b) \quad {\cI}^+_\alpha$ is dense open
decreasing with $\alpha$
\sn
\item[${{}}$]   $(c) \quad$ if $p \in \bbP_{{\cA}_*[\bold V]}$
\then \, for every large enough $\alpha < \lambda,p \notin {\cI}^+_\alpha$
\sn
\item[${{}}$]   $(d) \quad$ if $p \in \bbP_{\cA_*[\bold V]}$ and
  $\alpha < \lambda$ \then \, there is $q \in \cI_\alpha$ such that

\hskip25pt  $\bbP_{\cA_*[\bold V]} \models ``p \le q"$.
\end{enumerate}
\mn
Also clearly the sequence $\langle \cI_\alpha:\alpha <
\lambda\rangle$ belongs to $N_{\langle \rangle}$ hence 
 if $\alpha \in \lambda \cap N_{\name \eta[\bold G]}$ then
${\cI}_\alpha \in N_{\name \eta[\bold G]}$ and the set $\{p \in
{\cT} \cap N_{\name \eta[\bold G]}:p \le_{\cT} p_*$ and $p \in
\cT_\alpha\}$ is not empty.

Now
\mn
\begin{enumerate}
\item[$(*)_2$]   in $\bold V[\bold G]$ the following functions $h_\bullet,h_*$ 
are well defined
\sn
\begin{enumerate}
\item[$(a)$]   $\Dom(p_\bullet) = \Dom(h_*) = N_{<>} \cap {\gh}$
\sn
\item[$(b)$]  $h_\bullet(\gamma)$ is the unique $p \in 
N_{\name\eta[\bold G]} \cap \cT$ of level 
$\gamma$ which is $\le_{\cT} p_*$
\sn
\item[$(c)$]  if $\gamma < \gh$ then 
$h_*(\gamma) = h_{\gamma +1}(h_\bullet(\gamma +1))$
\end{enumerate}
\sn
\item[$(*)_3$]   if $\alpha \in \gh \cap N_{\name \eta[\bold G]}$
then $h_*(\alpha) \in N_{\name \eta[\bold G]} \cap {\gh} =
N_{<>} \cap {\gh}$
\end{enumerate}
\mn
also by the choice of $\bar h$ (and genericity) clearly
\mn
\begin{enumerate}
\item[$(*)_4$]   Rang$(h_*)$ is equal to $u := (2^{\aleph_0}) \cap 
N_{\name \eta[\bold G]}$.
\end{enumerate}
\mn
Lastly,
\mn
\begin{enumerate}
\item[$(*)_5$]   $h_* \in \bold V$.
\end{enumerate}
\mn
[Why?  As its domain, $N_{<>} \cap {\gh}$ belongs to $\bold V$ and
$h_*(\gamma)$ is defined from $\langle \cT,\bar h,\gamma,p_*\rangle \in
\bold V$ and $\cT$ is a tree.]
\mn
\begin{enumerate}
\item[$(*)_6$]   $(a) \quad$ from $u := \lambda \cap 
N_{\name \eta[\bold G]}$ we can define $\name \eta[\bold G]$
\sn
\item[${{}}$]   $(b) \quad u = \cup\{N_{\name \eta \rest n[\bold
G]} \cap \lambda:n < \omega\}$.
\end{enumerate}
\mn
[Why?  By the choice of $\bar N$.]

Together we get that $\name \eta[\bold G] \in \bold V$, contradiction.
\end{PROOF}

\begin{claim}
\label{2b.23}
We have $\neg \text{\rm Pr}_1(\bbQ,\bbP_{\cA_*[\bold V]})$ \when \,
\mn
\begin{enumerate}
\item[$(a)$]  $2^{\aleph_0} \ge \lambda = \text{\rm cf}(\lambda) >
\kappa = \gh$
\sn
\item[$(b)$]  $\alpha < \lambda \Rightarrow 
\cf([\alpha]^{\le \kappa},\subseteq) < \lambda$
\sn
\item[$(c)$]  $\bbQ$ is $(\lambda,\kappa)$-newly proper.
\end{enumerate}
\end{claim}

\begin{PROOF}{\ref{2b.23}}
Similar to \ref{2b.10}.
\end{PROOF}

\begin{conclusion}
\label{2b.26}
If $\gh < 2^{\aleph_0}$ and $\bbQ$ is a $(\gh^+,\gh)$-newly proper
\then \, $\neg\text{\rm Pr}_1(\bbQ,\bbP_{\cA_*[\bold V]})$.
\end{conclusion}
\newpage

\section {General sufficient conditions} \label{general}

\begin{claim}
\label{3c.1}  
Assume $\bold V \models \CH$.

If $\bbQ$ is c.c.c. \then \, {\rm Pr}$_2(\bbQ,\bbP_{{\cA}_*[\bold V]})$.
\end{claim}

\begin{remark}
1) This works replacing $\bbP_{{\cA}_*[\bold V]}$ by any
$\aleph_1$-complete $\bbP$ and strengthening the conclusions to
Pr$_1$, see \ref{3c.3}.

\noindent
2) See Definition \ref{1a.7}(1).
\end{remark}

\begin{PROOF}{\ref{3c.1}}
Let $\bbP = \bbP_{\cA_*[\bold V]}$.  Clearly it suffices to prove:
\mn
\begin{enumerate}
\item[$(*)$]   if $r \in \bbP$ and $\Vdash_{\bbQ} ``\name{\cI}$ 
is a dense open subset of $\bbP$" \then \, there is $r'$ such that:
\sn
\begin{enumerate}
\item[$(a)$]   $r \le_{\bbP} r'$
\sn
\item[$(b)$]   $\Vdash_{\bbQ} ``r' \in \name{\cI} \subseteq \bbP"$.
\end{enumerate}
\end{enumerate}
\mn
Why $(*)$ holds?  We try (all in $\bold V$) to choose
$(r_\alpha,q_\alpha)$ by induction on $\alpha < \omega_1$ but choosing
$q_\alpha$ together with $r_{\alpha +1}$ such that:
\mn
\begin{enumerate}
\item[$\circledast$]   $(a) \quad r_0 = r$
\sn
\item[${{}}$]   $(b) \quad r_\alpha \in \bbP$ is $\le_{\bbP}$-increasing
\sn
\item[${{}}$]   $(c) \quad q_\alpha \in \bbQ$
\sn
\item[${{}}$]   $(d) \quad q_\alpha,q_\beta$ are incompatible in
$\bbQ$ for $\beta < \alpha$
\sn
\item[${{}}$]   $(e) \quad q_\alpha \Vdash_{\bbQ} ``r_{\alpha +1} 
\in \name{\cI}"$.
\end{enumerate}
\mn
We cannot succeed in carrying the induction $\omega_1$ many steps
because $\bbQ \models$ c.c.c.

For $\alpha = 0$ no problem as only clause (a) is relevant.

For $\alpha$ limit - easy as $\bbP$ is $\aleph_1$-complete (and the only
relevant clause is (b)).

For $\alpha = \beta +1$, we first ask:
\bigskip

\noindent
\underline{Question}:  Is $\langle q_\gamma:\gamma < \beta\rangle$ a maximal
antichain of $\bbQ$?
\medskip

\noindent
\underline{If yes}, then $r_\beta$ is as required in $(*)$ on $r'$;
why?  if $\bold G_{\bbQ} \subseteq \bbQ$ is generic over $\bold V$ to 
which $r_\beta$ belongs, \then \, for some $\gamma <
\beta,q_\gamma \in \bold G_{\bbQ}$ hence $r_{\gamma +1} \in 
\name{\cI}[\bold G_Q]$ but $\name{\cI}[\bold G_{\bbQ}]$ is a
dense subset of $\bbP$ and is \underline{open} and $r_{\gamma +1}
\le_{\bbP} r_\beta$ so $r_\beta \in \name{\cI}[\bold G_{\bbQ}]$.

\noindent
\underline{If no}, let $q^\beta \in \bbQ$ be incompatible with $q_\gamma$ for
every $\gamma < \beta$.  Recalling $\Vdash_{\bbQ} ``\name{\cI}$ is
\underline{dense} and open" the set $X_\beta = \{r \in \bbP$: for some
$q,q^\beta \le_{\bbQ} q$ and $q \Vdash ``r \in \name{\cI}"\}$ is a
dense subset of $\bbP$ hence there is a member of $X_\beta$ above
$r_\beta$, let $r_\alpha$ be such member.  By $r_\alpha \in X_\beta$,
there is $q,q^\beta \le q$ such that $q \Vdash ``r_\alpha \in
\name{\cI}"$.  
So we choose $q_\beta$ as such $q$, so we can carry the induction
step.  

As said above we cannot carry the induction for all $\alpha <
\omega_1$ because then $\{q_\alpha:\alpha < \omega_1\}$ contradicts
``$\bbQ$ satisfies the c.c.c."  So for some $\alpha$ we cannot
continue, $\alpha$ is neither 0 nor limit hence for some $\beta,\alpha 
= \beta +1$.  So the answer to the question is yes, hence we
get the desired conclusion of $(*)$.
\end{PROOF}

We can weaken the demand on the second forcing 
(above, it is $\bbP_{{\cA}_*[\bold V]})$.
\begin{claim}
\label{3c.3}  
If (A) then (B) where:
\mn
\begin{enumerate}
\item[(A)]
\begin{enumerate}
\item[(a)]  $\bbP,\bbQ$ are forcing notions
\sn
\item[(b)]  $\bbQ$ is c.c.c. moreover $\Vdash_{\bbP} ``\bbQ$ is c.c.c."
\sn
\item[(c)]   forcing with $\bbP$ adds no new
$\omega$-sequences,\footnote{if you assume $\bbP$ is proper,
$\lambda = \aleph_0$ the proof may be easier to read} from $\lambda$
\sn
\item[(d)]  $\bbQ$ has cardinality $\le \lambda$ 
\end{enumerate}
\sn
\item[(B)]
\begin{enumerate}
\item[(a)]  if $\bbP$ is proper in $\bold V$ \then \, 
{\rm Pr}$_2(\bbQ,\bbP)$ 
\sn
\item[(b)]  for every $\bbQ$-name $\name {\cI}$ of a
dense open subset of $\bbP$, the set $\cJ$ is dense and open in $\bbP$ where:
\sn
\begin{enumerate}
\item[$(*)$]  $\cJ = \cJ_{\name{\cI}}$ is the set of $r \in \bbP$ such
  that some $\bar q$ witnesses it, i.e. witness it belongs to $\cJ$
  which means:
\sn
\begin{itemize}
\item  $\bar q = \langle q_\alpha:\alpha < \alpha_*\rangle$ is a
  maximal antichain of $\bbQ$
\sn
\item  for each $\alpha < \alpha_*$, the set $\{r' \in \bbP:q_\alpha
 \Vdash ``r' \in \name{\cI}"\}$ is an open subset of $\bbP$ dense above $r$.
\end{itemize}
\end{enumerate}
\end{enumerate}
\end{enumerate}
\end{claim}

\begin{PROOF}{\ref{3c.3}}
First, we prove clause (b); so fix $\name{\cI}$ and $\cJ$ as there.
Let $\langle q_\varepsilon:\varepsilon < \kappa := |\bbQ|\rangle$ 
list $\bbQ$.

For every $r \in \bbP$ we define a sequence $\eta_r$ of ordinals $<
\kappa \le \lambda$ as follows:
\mn
\begin{enumerate}
\item[$\circledast_1$]   $\eta_r(\alpha)$ is the minimal ordinal
$\varepsilon < \kappa$ such that (so $\ell g(\eta_r) = \alpha$ when
  there is no such $\varepsilon$):
\sn
\begin{enumerate}
\item[$(a)$]   $q_\varepsilon \Vdash ``r \in \name{\cI}"$
\sn
\item[$(b)$]   if $\beta < \alpha$ then
$q_\varepsilon,q_{\eta_r(\beta)}$ are incompatible in $\bbQ$.
\end{enumerate}
\end{enumerate}
\mn
Now
\mn
\begin{enumerate}
\item[$\circledast_2$]   $(a) \quad \eta_r$ is well defined
\sn
\item[${{}}$]   $(b) \quad \ell g(\eta_r) < \omega_1$.
\end{enumerate}
\mn
[Why?  Obviously $\eta_r$ is a well defined sequence of ordinals,
i.e. clause (a) and clause (b) holds because $\bbQ \models$ c.c.c.]

Note
\mn
\begin{enumerate}
\item[$\circledast_3$]   if $r_1 \le_{\bbP} r_2$ then 
\underline{either} $\eta_{r_1}
\trianglelefteq \eta_{r_2}$ \underline{or} for some $\alpha < \ell
g(\eta_{r_1})$ we have

\[
\eta_{r_1} \rest \alpha = \eta_{r_2} \rest \alpha
\]

\[
\eta_{r_1}(\alpha) > \eta_{r_2}(\alpha).
\]
\end{enumerate}
\mn
[Why?  Think about the definition.]

For $s \in \bbP$ let $\eta'_s$ be $\cap\{\eta_{s_1}:s \le_{\bbP}
s_1\}$, i.e. the longest common initial segment of $\{\eta_{s_1}:s
\le_{\bbP} s_1\}$; clearly $s_1 \le_{\bbP} s_2 \Rightarrow \eta'_{s_1}
\trianglelefteq \eta'_{s_2}$.  So 
\mn
\begin{enumerate}
\item[$\circledast_4$]   $\name \eta^* = \cup\{\eta'_s:s \in 
\name {\bold G}_{\bbP}\}$ is a $\bbP$-name of a sequence of ordinals
$<\kappa$ such that $\langle q_{\name\eta^*(i)}:i < \ell
g(\name\eta^*)\rangle$ is a sequence of pairwise incompatible members  
of $\bbQ$.
\end{enumerate}
\mn
But by clause (A)(b) of the claim, forcing with 
$\bbP$ preserve ``$\bbQ \models$ c.c.c.", 
so $\ell g(\name \eta^*)$ is countable in
$\bold V[\bold G_{\bbP}]$.  By clause (A)(c) of the claim,
forcing by $\bbP$ adds no new $\omega$-sequences to
$\kappa = |\bbQ|$ (and $\bbQ$ is infinite) and $\bold V[\bold G_{\bbP}]$ has 
the same $\aleph_1$ as $\bold V$, so
\mn
\begin{enumerate}
\item[$\circledast_5$]   $\name \eta^*$ is a
sequence of countable length of ordinals $< \kappa$ so is old.
\end{enumerate}
\mn
Hence
\mn
\begin{enumerate}
\item[$\circledast_6$]  the following set is dense open in $\bbP$
\[
{\cJ} = \{r \in \bbP:r \text{ forces in } \bbP \text{ that } 
\name \eta^* = \eta^*_r \text{ for some } \eta^*_r \in \bold V\}
\]
\end{enumerate}
\mn
As for clause (a), let $\chi,N,q_1,r_1$ be as in the assumption
of $(*)_1$ of \ref{1a.3}, so $\bbP,\bbQ \in N$.  We have to find
$q_2,r_2$ as there.

Let $q_2 = q_1$ and let $r_2 \in \bbP$ be $(N,\bbP)$-generic and above
$r_1$, exists as $\bbP$ is a proper forcing in $\bold V$.

We shall show that $(r_2,q_2)$ is as required, i.e. $q_2 \Vdash_{\bbQ}
``r_2$ is $(N[\name{\bold G}_{\bbQ}],\bbP)$-generic".  Let $\bold
G_{\bbQ} \subseteq \bbQ$ be generic over $\bold V$ such that $q_2 \in
\bold G_{\bbQ}$ and we should prove that $\bold V[\bold G_{\bbQ}]
\models `` r_2$ is $(N[\bold G_{\bbQ}],\bbP)$-generic".  So let $\cI
\in N[\bold G_{\bbQ}]$ be a dense open subset of $\bbP$, and we should
prove that $\bold V[\bold G_{\bbQ}] \models ``\cI \cap N[\bold
G_{\bbQ}]$ is pre-dense above $r_2$". 

It suffices to prove:
\mn
\begin{enumerate}
\item[$(*)$]  if $r_2 \le_{\bbP} r_3$ then $r_3$ is compatible (in
    $\bbP$) with some $r \in \cJ \cap N$.
\end{enumerate}
\mn
So fix $r_3 \in \bbP$; by the definition of $N[\bold G_{\bbQ}]$ there
is a $\bbQ$-name $\name{\cI}$ such that $\cI = \name{\cI}[\bold
G_{\bbQ}]$, for some $\name{\cI} \in N$;  \wilog \, $\Vdash_{\bbQ}
``\cI$ is a dense open subset of $\bbP$".  Let $\cJ = \cJ_{\name{\cJ}}
= \{r \in \bbP:r$ has an $\name{\cI}$-witness $\bar q_* = \langle
q^*_\alpha:\alpha < \alpha_*\rangle\}$, see clause (B)(b) of the claim.
Clearly $\cJ \in N$ hence $\cJ \cap N$ is pre-dense in $\bbP$ over
$r_2$ hence also over $r_3$ hence there are $r_4,r_5 \in \bbP$ such
that $r_3 \le_{\bbP} r_5,r_4 \le_{\bbP} r_5$ and $r_4 \in N \cap
\cJ$.  By the definition of $\cJ$ there is an $\name{\cI}$-witness
$\bar q_* = \langle q^*_\alpha:\alpha < \alpha_*\rangle$ for $r_4 \in
\cJ$.

But $\cI,r_4 \in N$ hence \wilog \, $\bar q_* \in N$ and $\bar q_*$
has countable length, so $\{q^*_\alpha:\alpha < \alpha_*\} \subseteq
N$.  As $\bar q_*$ is a witness, necesarily it is a maximal antichain
of $\bbQ$ hence for some $\alpha < \alpha_*$ we have $q^*_\alpha \in
\bold G_{\bbQ}$, as $\bar q_*$ is a witness for $r_4 \in
\cJ_{\name{\cI}}$, necessarily $\cI_1 = \{r \in \bbP:q^*_\alpha
\Vdash_{\bbQ} ``r \in \name{\cI}"\}$ is an open subset of $\bbP$ dense
above $r_4$.

Clearly $\cI_1 \in N$ is an open subset of $\bbP$, dense above $r_4$
and $r_4 \le_{\bbP} r_5$ hence $\cI_1 \cap N$ is pre-dense above $r_5$
hence there are $r_6 \le_{\bbP} r_7$ from $\bbP$ such that $r_6 \in
\cI_1 \cap N$ and $r_5 \le_{\bbP} r_7$.

Clearly $r_6 \in \cI[\bold G_{\bbQ}] \cap N$ and 
$r_6$ is compatible with $r_3$ in
$\bbP$, so we are done proving $r_2$ is $(N[\bold G_{\bbQ}],\bbP)$-generic. 

So we are done.
\end{PROOF}

\begin{remark}
\label{3c.5}
In \ref{3c.1}, \ref{3c.3} we can replace ``c.c.c." by ``strongly proper".

\underline{But} such $\bbQ$ preserves ``$({}^\omega 2)^{\bold V}$-non-meagre".
\end{remark}

\begin{claim}
\label{3c.7}  
1) There is a proper forcing $\bbQ$ which forces ``${\bbP}_{\cA_*}
[\bold V]$ as a forcing notion is not proper", 
(i.e. $\neg\text{\rm Pr}_1(\bbQ,\bbP))$.

\noindent
2) Even (A) of \ref{1a.10}(3) fails, i.e. 
$\neg\text{\rm Pr}_5(\bbQ,\bbP_{\cA_*}[\bold V])$. 
\end{claim}

\begin{PROOF}{\ref{3c.7}}
We use the proof of \cite[Ch.17,Sec.2]{Sh:f} and see 
references there.  We repeat in short.

We use a finite iteration so let $\bbP_0$ be the trivial forcing
notion, $\bbP_{k+1} = \bbP_k * \name{\bbQ}_k$ for $k \le 3$ and the
$\bbP_k$-name $\name{\bbQ}_k$ is defined below.
\bigskip

\noindent
\underline{Step A}:  $\bbQ_0 = \text{ Levy}(\aleph_1,2^{\aleph_0})$ so
$\Vdash_{\bbQ_0}$ ``CH".
\bigskip

\noindent
\underline{Step B}:  $\bbQ_1$ is Cohen forcing.
\bigskip

\noindent
\underline{Step C}:  In $\bold V^{\bbP_2},\bbQ_2$ in the Levy collapse of
$2^{2^{\aleph_0}}$ to $\aleph_1$, i.e. $\bbQ_2 = 
\text{ Levy}(\aleph_1,\beth_2)^{\bold V[\bbP_2]}$.
\bigskip

\noindent
\underline{Step D}: Let ${\cT} = ({}^{(\omega_1 >)}\omega_1)
^{\bold V[\bbP_1]} = ({}^{(\omega_1 >)} \omega_1)^{\bold V[\bbP_0]}$ 
be a tree, so we know that $\lim_{\omega_1}(\cT)^{\bold V[\bbP_1]} = 
\lim_{\omega_1}(\cT)^{\bold V[\bbP_2]} = 
\lim_{\omega_1}(\cT)^{\bold V[\bbP_3]}$ hence has cardinality
$\aleph_1$ in $\bold V^{\bbP_3}$ and
\mn
\begin{enumerate}
\item[$(*)_1$]  in $\bold V^{\bbP_1},\cT$ is isomorphic to a dense
subset of $\bbP_{\cA_*[\bbP_1]} = \bbP_{\cA_*[\bbP_0]}$.
\end{enumerate}
\mn
So in $\bold V^{\bbP_3}$ there is a list $\langle
\eta^*_\varepsilon:\varepsilon < \omega_1\rangle$ of 
$\lim_{\omega_1}(\cT)^{\bold V[\bbP_1]}$ and let 
$\langle \eta^*_\varepsilon \rest
[\gamma_\varepsilon,\omega_1):\varepsilon < \omega_1\rangle$ be
pairwise disjoint end segments so $\gamma_\varp < \omega_1,\langle
\gamma_\varp:\varp < \omega_1\rangle \in \bold V^{\bbP_3}$ and
$\varp_1 < \varp_2 < \omega_1 \wedge \beta_1 \in
[\gamma_{\varp_1},\omega_1) \wedge \beta_2 \in
[\gamma_{\varp_2},\omega_1) \Rightarrow \eta^*_{\varp_1} \rest
\gamma_1 \ne \eta^*_{\varp_2} \rest \gamma_2$.
\bigskip

\noindent
\underline{Step E}:  In $\bold V^{\bbP_3}$ there is $\bbQ_3$, 
a c.c.c. forcing notion specializing ${\cT}$ in the sense of
\cite{Sh:74}, i.e. there is $h_* \in \bold V^{\bbP_4}$ such that
 $h_*:{\cT} \rightarrow \omega,h_*$ is increasing
in ${\cT}$ except being constant on each end segment 
$\name\eta^*_\varepsilon \rest
[\gamma_\varepsilon,\omega_1)$ for $\varepsilon < \omega_1$, 
i.e. $\rho <_{\cT} \nu \wedge h_*(\rho) = h_*(\nu) 
\Rightarrow (\exists \varepsilon)[\rho,\nu \in
\{\eta^*_\varepsilon \rest \gamma:\gamma \in
[\gamma_\varepsilon,\omega_1)\}$.

Now
\mn
\begin{enumerate}
\item[$\boxtimes$]   after forcing with $\bbP_4 = \bbQ_0 * \bbQ_1 *
\name{\bbQ}_2 * \name{\bbQ}_3$, i.e. in $\bold V^{\bbP_4}$ 
the forcing notion $\bbP_{\cA_*[\bold V]}$ is not proper, 
in fact it collapses $\aleph_1$.
\end{enumerate}
\mn
Why?  Recall $(*)_1$ and note
\mn
\begin{enumerate}
\item[$(*)_2$]   ${\cI}_n := \{\rho \in {\cT}:(\forall \nu)(\rho
\le_{\cT} \nu \rightarrow h_*(\nu) \ne n\}$ is dense open in ${\cT}$
\end{enumerate}
\mn
and trivially
\mn
\begin{enumerate} 
\item[$(*)_3$]    $\bigcap\limits_{n} {\cI}_n = \emptyset$; in fact if
$\bold G \subseteq {\cT}$ is generic, \then \,:
\sn
\begin{enumerate} 
\item[$(A)$]    $\bold G$ is a branch of ${\cT}$ of order type
$\omega^{\bold V}_1$ let its name be $\langle \name\rho_\gamma:
\gamma < \omega_1\rangle$
\sn
\item[$(B)$]    letting $\name \gamma_n = \text{ Min}\{\gamma <
\omega_1:\name \rho_\gamma \in {\cI}_n\}$ we have $\Vdash_{\cT}
``\{\name \gamma_n:n < \omega\}$ is unbounded in $\omega_1$".
\end{enumerate}
\end{enumerate}
\end{PROOF}
\newpage

\bibliographystyle{alphacolon}
\bibliography{lista,listb,listx,listf,liste,listz}

\end{document}